\documentclass[article, lineno]{biometrika}

\usepackage{amsmath}

\usepackage{times}
\usepackage{bm}
\usepackage{natbib}

\usepackage[plain,noend]{algorithm2e}

\makeatletter
\renewcommand{\algocf@captiontext}[2]{#1\algocf@typo. \AlCapFnt{}#2} % text of caption
\def\@algocf@capt@plain{top}
\renewcommand{\algocf@makecaption}[2]{%
  \addtolength{\hsize}{\algomargin}%
  \sbox\@tempboxa{\algocf@captiontext{#1}{#2}}%
  \ifdim\wd\@tempboxa >\hsize%     % if caption is longer than a line
    \hskip .5\algomargin%
    \parbox[t]{\hsize}{\algocf@captiontext{#1}{#2}}% then caption is not centered
  \else%
    \global\@minipagefalse%
    \hbox to\hsize{\box\@tempboxa}% else caption is centered
  \fi%
  \addtolength{\hsize}{-\algomargin}%
}
\makeatother

\usepackage{graphicx,psfrag,epsf}
\usepackage{enumerate}
\usepackage{natbib}
\usepackage[caption=false]{subfig}
\usepackage{amsmath}
\usepackage{amssymb}
\usepackage{mathtools}
\usepackage{amsfonts}
\usepackage{mathptmx}
\usepackage
[acronym,nowarn,section,nonumberlist]{glossaries}
\glsdisablehyper{}

\usepackage{color} % pre-req for below macros
\usepackage{xcolor}
\usepackage{"./macro/GrandMacros"}
\usepackage{"./macro/Macro_BIO235"}

\usepackage{tocstyle}

\newcommand{\indep}{\rotatebox[origin=c]{90}{$\models$}}

\newacronym{UQ}{UQ}{uncertainty quantification}

\newacronym{MAP}{MAP}{maximum a posterior}

\newacronym{CDF}{CDF}{cumulative density function}
\newacronym{PDF}{PDF}{probability density function}
\newacronym{MGF}{MGF}{moment generating function}

\newacronym{GLM}{GLM}{generalized linear model}
\newacronym{GAM}{GAM}{generalized additive model}
\newacronym{MLE}{MLE}{maximum likelihood estimation}
\newacronym{CRPS}{CRPS}{\textit{continous ranked probability score}}
\newacronym{RMSE}{RMSE}{root mean squared error}
\newacronym{MSE}{MSE}{mean squared error}

\newacronym{KL}{KL}{Kullback-Leibler}
\newacronym{CvM}{CvM}{Cram\'{e}r von Mises}
\newacronym{ELBO}{ELBO}{\emph{evidence lower bound}}

\newacronym{BNE}{BNE}{\textit{Bayesian Nonparametric Ensemble}}
\newacronym{CVI}{Calibrated VI}{\textbf{Calibrated Variational Inference}}

\newacronym{RKHS}{RKHS}{reproducing kernel Hilbert space}
\newacronym{RBF}{RBF}{radial basis function}
\newacronym{CI}{CI}{Coverage Index}

\newacronym{GP}{GP}{Gaussian process}
\newacronym{CGP}{CGP}{\textit{constrained Gaussian process}}

\begin{document}

%\jname{Biometrika}
%%% The year, volume, and number are determined on publication
%\jyear{2018}
%\jvol{103}
%\jnum{1}
%%% The \doi{...} and \accessdate commands are used by the production team
%%\doi{10.1093/biomet/asm023}
%\accessdate{Advance Access publication on 31 July 2018}
%
%%% These dates are usually set by the production team
%\received{2 January 2017}
%\revised{1 April 2017}

%% The left and right page headers are defined here:
\markboth{}{}

%% Here are the title, author names and addresses
\title{
Gaussian Process Regression and Classification under Mathematical Constraints with Learning Guarantees
%Constrained Gaussian Process for Regression and Classification with Statistical Guarantees
}

\author{Jeremiah Z. Liu* \thanks{* Work in progress.}}
\affil{Department of Biostatistics, Harvard University,\\ Boston, MA, USA \email{zhl112@mail.harvard.edu.}}

\maketitle

\begin{abstract}
We introduce \gls{CGP}, a Gaussian process model for random functions that allows easy placement of mathematical constrains (e.g., non-negativity, monotonicity, etc) on its sample functions. %In Section \ref{sec:model}, we use \gls{CGP} to model the \gls{BNE}'s \gls{CDF} so the posterior respects the boundedness and monotonicity of a distribution function. 
\gls{CGP} comes with 
%an intuitive geometric interpretation as assigning probability mass to the intersection between convex sets, and has 
closed-form \gls{PDF}, and has the attractive feature that its posterior distributions for regression and classification are again \gls{CGP}s with closed-form expressions. Furthermore, we show that \gls{CGP} inherents the optimal theoretical properties of the Gaussian process, e.g. rates of posterior contraction, due to the fact that \gls{CGP} is an Gaussian process with a more efficient model space.
\end{abstract}


\newpage
%\section{Introduction}
%
%\subsection*{Notation}
%$\lesssim$\\
%$P$\\
%$\Hsc_k$\\
%$\phi$\\
%$\Phi$\\
%$\vee$\\

\section{The Idea}

Conceptually, the goal of a constrained Gaussian process (CGP) is to place mathematical constraints $\Csc$ on the sample paths $g$ of a Gaussian process, such that $g \in \Csc$. Operationally, $\Csc$ is defined through a system of $k$ linear constraints characterized by matrices $\bA_{k \times 1}$ and $\bb_{k \times 1}$, e.g.:
\begin{alignat*}{3}
    \mbox{sample} \quad & g \sim GP(0, k) \\
    \mbox{subject to} \quad & g \in \Csc 
    \qquad \qquad \mbox{where} \quad
    \Csc = \{g | \bA g + \bb \geq 0\}.
\end{alignat*}
For example, an inequality constraint $\Csc=\{g | g \geq 0 \}$ can be written as $\bA = 1, \bb = 0$, and a boundedness constraint $\Csc=\{g | 1 \geq g \geq 0 \}=\{g | g \geq 0,  1 - g \geq 0 \}$ can be written as $\bA = [1, -1]^\top, \bb = [0, 1]^\top$. Notice that since $g$ is a random variable, the constraint $\Csc$ is a random event with distribution $f(\Csc|g)$. Consequently, \gls{CGP} imposes constraints $\Csc$ on a Gaussian process by augmenting the original distribution of a \gls{GP} $f(g)$ with the additional likelihood term $P(\Csc|g)$.

As a result, given constraints $\Csc$ and a standard \gls{GP} with mean function $\mu$ and kernel function $k$, the likelihood function of \gls{CGP} is defined to be:
\begin{align}
    f(g|\Csc) \propto f(\Csc|g)f(g),
    \label{eq:cgp_lkhd}
\end{align}
where $f(g) = GP(g | \mu, k)$ is the likelihood for the standard \gls{GP}. $f(\Csc|g)$ is designed such that it assigns near zero probability to $g$'s that violates $\Csc$, and assign high probability otherwise. In this work, we use the Probit function (i.e. the \gls{CDF} of a standard Gaussian distribution) to represent the likelihood functions for each individual constraint in $\Csc$, and assume the $k$ constraints are independent conditional on $G$. As a result, $P(\Csc|G)$ is the \gls{CDF} of a standard multivariate Gaussian distribution:
\begin{align*}
    f(\Csc|g) \propto \Phi(\bA g + \bb | \bmu=\bzero, \bSigma=\bI)
\end{align*}
where we denote $\Phi(.|\bmu, \bSigma)$ and $\phi(.|\bmu, \bSigma)$ to be the \gls{CDF} and \gls{PDF} of a Gaussian distribution with mean $\bmu$ and covariance $\bSigma$. To keep the notations uncluttered, we will drop the terms $\bmu$ and $\bSigma$ unless they are different from the default values $\bmu=\bzero$ and $\bSigma=\bI$. Consequently, conditional on the function input $\bz$, the  likelihood function of \gls{CGP} is:
\begin{align}
    f(g|\Csc) 
    &\propto f(\Csc|\bg) * f(g), \nonumber \\
    &\propto \Phi(\bA g + \bb) * \phi(g |\mu,  k),
    \label{eq:cgp_cond_lkhd}
\end{align}
where $\bu$ and $\bK$ are the mean vector and the kernel matrix of the standard GP $f(g) = GP(g | \mu, k)$ evaluated at  $\bz$.
%\jzl{add geometric interpretation}

%\section{The Constrained Gaussian Process}
%\subsection{Definition}
%\jzl{add a rigorous definition here.}

\subsection{\gls{PDF}, \gls{MGF} and Moments}
Given constraints $\Csc = \{ g | \bA g + \bb \geq 0\}$ and a standard Gaussian process $GP(\mu, k)$, denote the Constrained Gaussian Process (CGP) as $CGP(\mu, k, \Csc)$. For $g \sim CGP(\Csc, \mu, k)$, given input $\bz$, the \gls{PDF} of $\bg = g(\bz)$ is:
\begin{align}
    f(\bg |\bu, \bK, \Csc) = 
    \frac{\phi(\bg |\bu,  \bK) * \Phi(\bA \bg + \bb)}{C}
    \quad \mbox{where} \quad 
    C = \Phi(\bA\bu + \bb|\Sigma = \bI + \bA\bK\bA^\top)
    \label{eq:cpg_pdf}
\end{align}
and the \gls{MGF} of $CGP(\Csc, \mu, k)$ is:
\begin{align}
    M_{\bg}(\bt) 
    %&= E(exp(\bt^\top\bg)|\Csc, \bu, \bK)  \nonumber\\
    &= exp(\bu^\top\bt + \frac{1}{2} \bt^\top \bK \bt ) * 
    \frac{
    \Phi(\bA\bu + \bb + \bA\bK\bt |\Sigma = \bI + \bA\bK\bA^\top)}{
    \Phi(\bA\bu + \bb|\Sigma = \bI + \bA\bK\bA^\top)}
    \label{eq:cpg_mgf}
\end{align}
where recall $\Phi(.|\bmu, \bSigma)$ and $\phi(.|\bmu, \bSigma)$ are the \gls{CDF} and \gls{PDF} of a Gaussian distribution with mean $\bmu$ and covariance $\bSigma$.

%\jzl{see if we can derive a even more interpretable form.}
As a result, we can derive the moments of $CGP(\Csc, \mu, k)$ by taking derivatives with respect to (\ref{eq:cpg_mgf}). In particular, the expression for the mean of CGP is
\begin{align}
    E(\bg|\bu, \bK, \Csc) 
    &=
    \bu + 
    \bK\bA^\top (\bI + \bA\bK\bA^\top)^{-1}
    E(\bb'| \bb' \geq -\bb, \Sigma=\bI + \bA\bK\bA^\top)
    \label{eq:cgp_mean}
\end{align}
In particular, the second term in (\ref{eq:cgp_mean}) corrects the original mean vector $\bmu$ with respect to the constraints such that $\bA E(\bg |\bu, \bK, \Csc) + \bb \geq 0$.
\begin{proof} 
See Section \ref{sec:cgp_pdf_mgf_proof}
\end{proof}

\section{Regression and Classification}
An important feature of \gls{CGP} is that the posterior distributions of $\bg$ in regression and classification have closed forms. Specifically, assuming $g \sim CGP(\mu, k, \Csc)$, for a regression model $y_i \stackrel{indep}{\sim} Normal(g(\bz_i), \sigma^2)$, the posterior distribution of $\bg$ is again a constrained Gaussian process $CGP(\mu', k', \Csc')$ with parameters:
\begin{align*}
%    f(\bg |\Csc, \bu, \bK, \by) &= 
%    \frac{\phi(\bg |\bmu',  \Sigma') * \Phi(\bA \bg + \bb)}{
%    \Phi(\bA\bmu' + \bb|\Sigma = \bI + \bA\bSigma'\bA^\top)}
%    \qquad \mbox{where} \quad 
%    \label{eq:cpg_pdf_posterior}
%    \\
    \bmu' &= \bu + \bK(\bK + \sigma^2\bI)^{-1}(\by - \bu), 
    \quad 
    \bK' = \bK - \bK(\bK + \sigma^2\bI)^{-1}\bK,
    \quad 
    \Csc' = \Csc,
\end{align*}
and the predictive distribution for a new observation $\bz^*$ also follows a constrained Gaussian process $CGP(\mu^*, k^*, \Csc^*)$ with parameters:
\begin{align*}
\bmu^* &= u(\bz^*) + k(\bz^*, \bz)(\bK + \sigma^2\bI)^{-1}(\by - \bu),\\
\bK^* &= k(\bz^*, \bz^*) -  k(\bz^*, \bz)(\bK + \sigma^2\bI)^{-1}k(\bz, \bz^*),\\
\Csc^* &= \Csc.
\end{align*}

For a classification model under the Probit link function $y_i \stackrel{indep}{\sim} Bernoulli(p_i)$ where $p_i = \Phi(g(\bz_i))$, the posterior distribution of $\bg$ is also a constrained Gaussian process $CGP(\mu', k', \Csc')$ with modified constraint $\Csc'$:
\begin{align*}
    \bmu' = \bu, 
    \qquad 
    \bK' = \bK, 
    \qquad 
    \Csc' = \Bigg\{ \bg \Big|
    \begin{bmatrix} \bA \\ \bD \end{bmatrix} \bg + \begin{bmatrix} \bb \\ \bzero \end{bmatrix} \geq 0
    \Bigg\}
\end{align*}
where $\bD = diag(2 * \by - 1)$ is a diagonal matrix of $\pm 1$'s corresponding to the observations. Notice that the $\bD\bg$ term in $\Phi$ can be interpreted as imposing addition data-based constraints to $g(\bz_i)$, i.e.  $g(\bz_i) \geq 0$ when $y_i = 1$ and $g(\bz_i) \leq 0$ otherwise, such that during estimation,  $g(\bz_i)$ is pushed toward positive/negative values depending on the value of the observation $y_i$. The predictive distribution for $g$ at new location $\bz^*$ is also a constrained GP $CGP(\mu^*, k^*, \Csc^*)$, with parameters:
\begin{align*}
\bmu^* &= u(\bz^*) + k(\bz^*, \bz)\bK^{-1}(\by - \bu),\\
\bK^* &= k(\bz^*, \bz^*) -  k(\bz^*, \bz)\bK^{-1}k(\bz, \bz^*),\\
\Csc^* &= \Bigg\{ \bg \Big|
        \begin{bmatrix} \bA \\ \bD \end{bmatrix} \bg + 
        \begin{bmatrix} \bb \\ \bzero \end{bmatrix} \geq 0
        \Bigg\}
\end{align*}
\begin{proof}
See Section \ref{sec:cgp_pdf_posterior_proof}
\end{proof}

\section{Posterior Concentration}

\subsection{Feasibility Condition}
For a target function $g^* \in \Hsc_k$ that is \textit{feasible} with respect to  constraint $\Csc$ (i.e. $g^* \in \Csc$), the interest of \gls{CGP} is to better estimate $g^*$ by shifting its probability mass toward the region where $f(\Csc|g)$ is high. Consequently, to measure the convergence of \gls{CGP}'s posterior toward $g^*$, we need a notion about the ``degree" of feasibility of $g^*$ with respect to the probablistic constraint $f(\Csc|g)$ specified by \gls{CGP}, so that we can decide if our configuration of the \gls{CGP} is ``compatible enough" with $g^*$ to guarantee fast speed of posterior convergence. In this work, we establish such notion by considering how robust $f(\Csc|g^*)$ is under a small amount of random perturbation $s$ on $g^*$:

\begin{definition}[$\epsilon$-feasiblity]
\label{def:feasibility}

Denote $truncGP(0, k, S)$ the truncated Gaussian process with mean zero, covariance kernel $k$ and its support truncated to be within the set $S$. For a small positive constant $\epsilon > 0$, denote $S_\epsilon = \{s | ||s|| \leq \epsilon\}$ a set of perturbation noises with maximum magnitude $\epsilon$.

For a function $g^* \in \Hsc_k$, we say $g^*$ is \textbf{$\epsilon$-feasible} with respect to the probablistic constraint $f(\Csc|g)$ if for the ``$\epsilon$-perturbed" function $g^*_\epsilon truncGP(0, k, g^* + S_\epsilon)$ and the random noise $s_\epsilon \sim truncGP(0, k, S_\epsilon)$, we have 
\begin{align*}
E \Big( f(\Csc|g^*_\epsilon) \Big) \geq E \Big( f(\Csc|s_\epsilon) \Big)
\end{align*}
i.e. the $\epsilon$-pertubed function $g^*_\epsilon$ is on average more feasible than the random noise $s_\epsilon$.
\end{definition}

The notion of $\epsilon$-feasiblity requires $g^*$ to be more feasible than the random noises under random perturbations of magnitude $\epsilon$. Notice that when $\epsilon$ is very small (as is the case for Theorem \ref{thm:cgp_conv}), the $\epsilon$-feasibility condition is essentially requiring $f(\Csc|g^*) \geq f(\Csc|\bzero)$, i.e. the target function $g^*$ should be no less feasible than the zero function $\bzero(x)=0$, which is the default mean function of a \gls{CGP} prior. 

\subsection{Posterior Concentration for General \gls{CGP} Prior}
For a target function $g^*$ that is $\epsilon$-feasible, we can show that the \gls{CGP} prior assigns sufficient probability mass around its neighborhood, which is important for guaranteeing reasonable speed of posterior convergence toward $g^*$ (see Lemma \ref{thm:feasibility_ineq} in Appendix). Furthremore, we can show that \gls{CGP} enjoys an theoretical guarantee in the posterior convergence toward a target functions that are $\epsilon$-feasible:

\begin{theorem}[Conditions for Posterior Consistency in \gls{CGP}]
Let $g$ be a Borel measurable, zero-mean constratined Gaussian random element in a separable Banach space $(\Bbb, ||.||)$ with \gls{RKHS} $(\Hsc_k, ||.||_{\Hsc_k})$. Define the concentration function $\psi_{g^*}(\epsilon) = \inf_{\hat{g} \in \Hsc_k, ||\hat{g} - g^*|| \leq \epsilon} ||\hat{g}||^2_{\Hsc_k} - log \, P(||g|| \leq \epsilon)$. 

For any number $\epsilon_n > 0$ satisfying $\psi_{g^*}(\epsilon_n) \leq n\epsilon_n^2$, and any constant $C \geq 1$ with $e^{-Cn\epsilon_n^2} < \frac{1}{2}*  \frac{E(f(\Csc|g))}{E(||g||_{\Hsc_k}^2)\vee 1}$, for $g^*$ a function contained in the closure of $\Hsc_k$ in $\Bbb$ that is $\epsilon_n$-feasible, there exists a measurable set $B_n \subset \Bbb$ such that:
\begin{align}
P(||g - g^*|| < 2\epsilon_n) &\geq e^{-n\epsilon_n^2} 
\tag{I}\label{eq:conv_cond_1} \\
P(g \not\in B_n) &\leq e^{-Cn\epsilon_n^2} 
\tag{II}\label{eq:conv_cond_2} \\
log \, N(2 \epsilon_n, B_n, ||.||) &\leq 2 C n \epsilon_n^2 
\tag{III}\label{eq:conv_cond_3} 
\end{align}
\label{thm:cgp_conv}
\end{theorem}
In above theorem, $B_n$ can be understood as the "large probability region" of a \gls{CGP} model, i.e., region where the posterior distribution will put  sufficiently large amount of probability mass in. Ideally, as the sample size $n$ grow, we hope this region to move quickly from the initial location to concentrate around the target function $g^*$. To this regard, the three conditions in Theorem \ref{thm:cgp_conv} describes how the \gls{CGP} prior behave with respect to the data-generating function $g^*$ and a  "model" $B_n$. Specifically, condition (\ref{eq:conv_cond_1}) requires the \gls{CGP} prior to put sufficient mass around $g^*$, condition (\ref{eq:conv_cond_2}) requires the prior to be not too big compared to the model $B_n$, and condition (\ref{eq:conv_cond_3}) puts a restriction on the size of the model $B_n$, in the sense that the size of $B_n$ when measured by the entropy number (i.e. the minimum number of balls of radius $3\epsilon_n$ needed to cover $B_n$) is upper bounded.

Similar to Theorem 2.1 of \cite{van_der_vaart_bayesian_2007} which outlines the general convergence conditions of Gaussian processes, the three conditions in Theorem \ref{thm:cgp_conv} can be matched one-to-one to the conditions for posterior convergence in \cite{ghosal_convergence_2000} (Theorem 2.1), with the exception that the distance measures in Theorem \ref{thm:cgp_conv} are defined to be $||.||$ (i.e. the norm for the Banach space $\Bbb$ with which the constrained Gaussian process is defined) rather than the typical  statistical distances (e.g. the Hellinger distance) that were used to measure convergence. Consequently, for a statistical problem under consideration (e.g. regression or classification), we can show convergence by showing that the statistical metric for measuring convergence is bounded by the Banach space norm $||.||$ and then invoke Theorem \ref{thm:cgp_conv}.

%\begin{corollary}[Posterior Concentration, Regression]
%
%\end{corollary}
%
%\begin{corollary}[Posterior Concentration, Classification]
%
%\end{corollary}
%
%\subsection{Concentration Rates for Concrete Kernel Families}
%
%
%
%\section{Application to Distribution Estimation}
%
%
%
%\section{Discussion}

\clearpage
\bibliographystyle{biometrika}
\bibliography{../../report}

\clearpage
\appendix
\appendix
\section{Derivation}
\subsection{Derivation for \gls{CGP}'s \gls{PDF}, \gls{MGF} and Mean}
\label{sec:cgp_pdf_mgf_proof}
To derive the expression of PDF in (\ref{eq:cpg_pdf}), only need to derive the constant term $C$ by integrating over \gls{CGP}'s unnormalized likelihood function in (\ref{eq:cgp_cond_lkhd}):
\begin{proof}
\begin{alignat}{3}
    C 
    &= \int_{\bg \in \real^n} \Phi(\bA \bg + \bb) * \phi(\bg |\bu,  \bK) d\bg 
    \nonumber\\
    &= P \big( \bZ \leq \bA \bg + \bb \big)  \qquad 
    &&\mbox{where}\quad 
    \bZ \sim \phi(0, \bI)
    \nonumber\\
    &= P \big( \bZ \leq \bA(\bGamma \bZ' + \bu) + \bb \big)  \qquad 
    &&\mbox{where}\quad 
    \bZ' \sim \phi(0, \bI), \bZ' \indep \bZ, \bA=\bGamma\bGamma^\top 
    \nonumber\\
    &= P \big( \bZ - \bA\bGamma \bZ' \leq \bA\bu + \bb \big)  \qquad 
    &&\mbox{notice}\quad 
    \bZ - \bA\bGamma \bZ' \sim \phi(0, \Sigma=\bI + \bA\bK\bA^\top)
    \nonumber\\
    &= \Phi(\bA\bu + \bb|\bSigma=\bI + \bA\bK\bA^\top)
    \label{eq:cdf_integral}
\end{alignat}
\end{proof}

Now derive the expression for \gls{MGF}:
\begin{proof}
\begin{alignat*}{3}
    M_{\bg}(\bt) 
    &= E(exp(\bt^\top \bg )|\Csc, \bu, \bK) \\
    &= \int_{\bg \in \real^n}
    exp(\bt^\top \bg ) * \phi(\bg |\bu,  \bK) * \Phi(\bA \bg + \bb)d\bg /C  \\
    &\propto \int_{\bg \in \real^n} 
    exp(- \frac{1}{2} \bg \bK^{-1} \bg + (\bK^{-1}\bu + \bt)^\top\bg ) * \Phi(\bA \bg + \bb)  d\bg /C 
    \qquad &&\mbox{complete square in exp term}
    \\
    &= \int_{\bg \in \real^n} 
    exp(\bu^\top\bt + \frac{1}{2} \bt^\top \bK \bt ) * 
    \phi(\bg|\mu = \bu + \bK\bt, \Sigma=\bK) * \Phi(\bA \bg + \bb)  d\bg /C
    \\
    &= exp(\bu^\top\bt + \frac{1}{2} \bt^\top \bK \bt ) * 
    \int_{\bg \in \real^n} 
    \phi(\bg|\mu = \bu + \bK\bt, \Sigma=\bK) * \Phi(\bA \bg + \bb) d\bg  /C
    \qquad &&\mbox{compute integration as in (\ref{eq:cdf_integral})}
    \\
    &= exp(\bu^\top\bt + \frac{1}{2} \bt^\top \bK \bt ) * 
    \Phi(\bA(\bu + \bK\bt) + \bb | \Sigma = \bI + \bA\bK\bA^\top)/C
    \qquad &&\mbox{plug in expression of $C$}
    \\
    &= exp(\bu^\top\bt + \frac{1}{2} \bt^\top \bK \bt ) * 
    \frac{
    \Phi(\bA\bu + \bb + \bA\bK\bt |\Sigma = \bI + \bA\bK\bA^\top)}{
    \Phi(\bA\bu + \bb|\Sigma = \bI + \bA\bK\bA^\top)}
\end{alignat*}
\end{proof}

Derive expression for \gls{CGP}'s mean in (\ref{eq:cgp_mean}):
\begin{proof}
\begin{align*}
    E(\bg|\bu, \bK, \Csc) 
    = \nabla_\bt M_\bg(\bt) \Big|_{\bt=\bzero}
    &=
    \bu - \bK\bA^\top
    \frac{
    \nabla \Phi(\bA\bu + \bb|\Sigma = \bI + \bA\bK\bA^\top)
    }{
    \Phi(\bA\bu + \bb|\Sigma = \bI + \bA\bK\bA^\top)
    }\\
    &=
    \bu - \bK
    \frac{
    \nabla_{\bu} \Phi(\bA\bu + \bb|\Sigma = \bI + \bA\bK\bA^\top)
    }{
    \Phi(\bA\bu + \bb|\Sigma = \bI + \bA\bK\bA^\top)
    }
\end{align*}
where $\nabla_\bu \Phi$ denotes the gradient of $\Phi$ with respect to $\bu$.

Furthermore, it can be shown that:
\begin{align*}
    \nabla_{\bu} \Phi(\bA\bu + \bb|\Sigma) 
    &=
    \int_{\bx \leq \bzero} \nabla_{\bu} \phi(\bx|\bmu = \bA\bu + \bb, \Sigma) d\bx \\
    &=
    \bA^\top \bSigma^{-1}
    \int_{\bx \leq \bzero} (\bx - \bA\bu) 
    \phi(\bx|\bmu = \bA\bu + \bb, \Sigma) d\bx
    \\
    &=
    \bA^\top \bSigma^{-1}
    \int_{\bx' \leq \bb} \bx' \phi(\bx'| \Sigma) d\bx,
\end{align*}
therefore 
\begin{align*}
    \frac{
    \nabla_{\bu} \Phi(\bA\bu + \bb|\Sigma = \bI + \bA\bK\bA^\top)
    }{
    \Phi(\bA\bu + \bb|\Sigma = \bI + \bA\bK\bA^\top)
    } &=
    \bA^\top \bSigma^{-1}
    \frac{
    \int_{\bx' \leq \bb} \bx' \phi(\bx'| \Sigma) d\bb'}{
    \int_{\bx' \leq \bb} \phi(\bx'| \Sigma) d\bb'    
    } \\
    &=  \bA^\top \bSigma^{-1} E(\bx'| \bx' \leq \bb, \Sigma)\\
    &=  -\bA^\top \bSigma^{-1} E(\bb'| \bb' \geq -\bb, \Sigma)
\end{align*}
Plugging above expression into the expression for $E(\bg|\bu, \bK, \Csc)$, we get
\begin{align*}
    E(\bg|\bu, \bK, \Csc) 
    &=
    \bu - \bK
    \frac{
    \nabla_{\bu} \Phi(\bA\bu + \bb|\Sigma = \bI + \bA\bK\bA^\top)
    }{
    \Phi(\bA\bu + \bb|\Sigma = \bI + \bA\bK\bA^\top)
    }\\
    &= 
    \bu +  \bK\bA^\top \bSigma^{-1} E(\bb'| \bb' \geq - \bb, \Sigma).
\end{align*}
This concludes the proof.
\end{proof}

\subsection{Derivation for \gls{CGP}'s posterior \gls{PDF}}
\label{sec:cgp_pdf_posterior_proof}
First derive the posterior distribution for regression model:
\begin{align*}
    y_i &\stackrel{indep}{\sim} Normal(g(\bz_i), \sigma^2) \\
    g &\sim CGP(\mu, k, \Csc)
\end{align*}

\begin{proof}
To derive the posterior \gls{PDF}, first write out its unnormalized form:
\begin{align*}
f(\bg|\Csc, \bu, \bK, \by) 
&\propto f(\by|\bg)f(\bg|\Csc, \bu, \bK)\\
&\propto exp(-\frac{1}{2\sigma^2}(\by - \bg)^\top(\by - \bg)) * \phi(\bg |\bu,  \bK) * \Phi(\bA \bg + \bb) \\
&\propto exp(-\frac{1}{2}\bg^\top(\bK^{-1} + \sigma^{-2}\bI)\bg + (\by - \bK^{-1}\bu)^\top\bg) *
\Phi(\bA \bg + \bb) \\
&\propto \phi(\bg|\mu', \Sigma')\Phi(\bA\bg + \bb)
\end{align*}
where $\mu'=(\bK^{-1} + \sigma^{-2}\bI)^{-1}(\frac{1}{\sigma^2}\by - \bK^{-1}\bu)$ and $\Sigma'=(\bK^{-1} + \sigma^{-2}\bI)^{-1}$. Using spectral decomposition of $\bK$, it can be shown easily that $\mu'=\bu + \bK(\bK + \sigma^2\bI)^{-1}(\by - \bu)$ and $\Sigma'=\bK - \bK(\bK + \sigma^2\bI)^{-1}\bK$ as is done in the standard \gls{GP} model \citep{rasmussen_gaussian_2006}. Consequently, the expression for the posterior \gls{PDF} is:
\begin{align*}
f(\bg|\Csc, \bu, \bK, \by) &= \frac{\phi(\bg|\mu', \Sigma')\Phi(\bA\bg + \bb)}{C}\\
C &= \int_{\bg \in \real^n} \phi(\bg|\mu', \Sigma')\Phi(\bA\bg + \bb) d\bg \\
&= \Phi(\bA\bmu' + \bb|\Sigma = \bI + \bA\bSigma'\bA^\top) \qquad \mbox{Using result from } (\ref{eq:cdf_integral})
\end{align*}
where recall $\mu'=\bu + \bK(\bK + \sigma^2\bI)^{-1}(\by - \bu)$ and $\Sigma'=\bK - \bK(\bK + \sigma^2\bI)^{-1}\bK$.
\end{proof}

Now derive the posterior distribution for classification model:
\begin{align*}
y_i &\stackrel{indep}{\sim} Bernoulli(p_i), p_i = \Phi(g(\bz_i)) \\
g &\sim CGP(\mu, k, \Csc)
\end{align*}

\begin{proof}
To derive the posterior \gls{PDF}, first write the expression for the likelihood $f(\by|\bg)$. Denote $0 \leq i^+ \leq n^+$ and $0 \leq i^- \leq n^-$ the index set corresponding to positive and negative observations, then:
\begin{align*}
f(\by|\bg) 
&\propto 
\prod_{i^+=1}^{n^+} \Phi(g(\bz_i)) * \prod_{i^-=1}^{n^-}(1-\Phi(g(\bz_j)))
\\
&=
\prod_{i^+=1}^{n^+} \Phi(g(\bz_i)) * \prod_{i^-=1}^{n^-}(\Phi(-g(\bz_j)))
\\
&= 
\prod_{i}^n \Phi  \big( (2 y_i - 1) * g(\bz_i) \big)
\\
&= \Phi(\bD \bg)
\end{align*}
where $\bD = diag(2 * \by - 1)$. Consequently, the posterior likelihood is:
\begin{align*}
f(\bg|\Csc, \bu, \bK, \by) 
&\propto f(\by|\bg)f(\bg|\Csc, \bu, \bK)\\
&\propto f(\bg|\Csc, \bu, \bK) \Phi(\bD \bg) \\
&\propto \phi(\bg|\mu', \Sigma')\Phi(\bA\bg + \bb) \Phi(\bD \bg) \\
&\propto \phi(\bg|\mu', \Sigma')\Phi(
\begin{bmatrix}\bA \\ \bD \end{bmatrix} \bg + \begin{bmatrix} \bb \\ \bzero \end{bmatrix}
)
\end{align*}
\end{proof}

\section{Additional Definitions}
\subsection{Full definition for $\epsilon$-feasibility}
\label{sec:cgp_feasibility}
We state the full definition of $\epsilon$-feasiblity:
\begin{definition}[$\epsilon$-feasiblity, full definition]
\label{def:feasibility_app}

Denote $truncGP(0, k, S)$ the truncated Gaussian process with mean zero, covariance kernel $k$ and its support truncated to be within the set $S$. For a small positive constant $\epsilon > 0$, denote $S_\epsilon = \{s | \; ||s|| \leq \epsilon\}$ a ``perturbation set" whose elements are random noise with maximum magnitude $\epsilon$. Also denote $\alpha_\epsilon = E\big( f(\Csc|s_\epsilon) \big)$ the ``average feasibility" of random noise for $s_\epsilon \sim truncGP(0, k, S_\epsilon)$.

For a function $g^* \in \Hsc_k$, we say $g^*$ is \textbf{$\epsilon$-feasible} with respect to the probablistic constraint $f(\Csc|g)$ if the ``$\epsilon$-perturbed" function $g^* + s\sim truncGP(g^*, k, g^* + S_\epsilon)$ satisfies any of the below conditions:
\begin{itemize}
\item ($\epsilon$-feasiblity, almost surely):
\begin{align*}
P\Big( f(\Csc|g^* + s) \geq \alpha_\epsilon \Big) = 1
\end{align*}
i.e. the $\epsilon$-perturbed function $g^* + s$ is almost always more feasible than the random noise $s_\epsilon$.
\item ($\epsilon$-feasibility, in probability):
\begin{align*}
P\Big( f(\Csc|g^* + s) < \alpha_\epsilon \Big) \leq \beta_\epsilon 
\end{align*}
i.e. the probability that the $\epsilon$-perturbed function $g^* + s$ is less feasible than random noise is upper bounded by a constant $\beta_\epsilon$. \\
In this work,  we set $\beta_\epsilon = 
P \Big( f(\Csc|g^* + s') < \alpha_\epsilon \Big| s' \in S_\epsilon \Big)$ for $g^* + s' \sim CGP(g^*, k, \Csc)$.
\item ($\epsilon$-feasibility, in expectation):
\begin{align*}
E \Big( f(\Csc|g^* + s) \Big) \geq \alpha_\epsilon 
\end{align*}
i.e. the $\epsilon$-perturbed function $g^* + s$ is on average more feasible than the random noise $s_\epsilon$.
\end{itemize}

\end{definition}

The definition of $\epsilon$-feasiblity requires $g^*$ to be ``feasible enough" for $\Csc$ such that it is  more feasible than the random noises even under random perturbation. The definition of "in probability" feasiblity relaxes the almost-sure feasibility by only requiring the probability of 
violating the strong feasibility is small, i.e., upper bounded by a value $\beta_\epsilon$. Notice that in the above definition, $\beta_\epsilon$  measures the probability of the feasiblity condition being violated by a ``specially perturbed"  function $g^* + s' \sim CGP(g^*, k, \Csc)$. Specifically, this function is similar to $g^*+s$ in that it perturbs $g^*$ using random noises with maximum magnitude $\epsilon$, but it is "special" in the sense that the noise distribution is specially designed such that $g^* + s'$ still respects the probablistic constraints $\Csc$ even after perturbation. Conseqeuntly, the weak feasibility condition essentially states that the perturbed $g^*$ should always be more feasible than the random noise regardless of how it is perturbed (i.e. if the perturbation noise respects $\Csc$ or not). Also notice that when $\epsilon$ is very small (as is the case for Theorem \ref{thm:cgp_conv}), the $\epsilon$-feasibility condition is essentially requiring $f(\Csc|g^*) \geq f(\Csc|\bzero)$, i.e. the target function $g^*$ needs to be as or more feasible than the prior mean function $\bzero(x)=0$. 

\section{Proof}
\subsection{Proof for \gls{CGP}'s posterior consistency}
\label{sec:cgp_convergence_proof}
%Recall the definition of $\epsilon$-feasiblity considers how robust $g^*$ is with respect to the probablistic constraint $f(\Csc|g)$ under a small amount of random perturbation $s$. Specifically, it requires $g^*$ to be "feasible enough" such that even under perturbation $\epsilon \in S_\epsilon$, the value of $f(\Csc|g^* + s)$ evaluated at the purturbed function $g^* + s$ is still  lower bounded by a value $\alpha_\epsilon$. In this work, we set this lower bound to be the expectation of $f(\Csc|s)$ evaluated at the random perturbation noise $s \in S_{\epsilon}$, i.e. the ``average feasibility" of some random noise $s$. Notice that for very small $\epsilon$ (as is the case for Theorem \ref{thm:cgp_conv}), the $\epsilon$-feasibility condition is essentially requiring $f(\Csc|g^*) \geq f(\Csc|\bzero)$, i.e. $g^*$ needs to be no less feasible than the zero function $\bzero(x)=x$. 

For a target function $g^*$ that is $\epsilon$-feasible (in any sense as defined in Definition \ref{def:feasibility_app}), we can show that the \gls{CGP} prior assigns sufficient probability mass around its neighborhood, which is important for guaranteeing reasonable speed of posterior convergence toward $g^*$:

\begin{lemma}[$\epsilon$-feasible function receives sufficient mass from \gls{CGP}]
Let $g$ be distributed as a zero-mean \gls{CGP} with covariance kernel $k$ and constraint $\Csc$. Consider a small ball $S_\epsilon =\{g \big| ||g|| \leq \epsilon \}$, then for a $\bepsilon$-feasible function $g^* \in \Hsc_k$, we have:
\begin{align*}
P(g: ||g - g^*|| \leq \epsilon) \geq 
exp \Big(-||g||^2_{\Hsc_k} \Big) * 
P(g: ||g|| \leq \epsilon)
\end{align*}
\label{thm:feasibility_ineq}
\end{lemma}
Intuitively, above result implies that for a function $g^*$ is $\epsilon$-feasible, the neighborhood surrounding $g^*$ with radius $\epsilon$ is always receive "sufficient" amount of probability mass from the \gls{CGP} prior, in the sense that the amount of mass received is bounded away from zero by a function of $P(S_\epsilon)$.

This result is analogous to the shifted-ball inequality for the standard Gaussian measures by \cite{kuelbs_gaussian_1994} (Theorem 2), which played a key role in establishing the posterior convergence conditions for Gaussian processes \citep{van_der_vaart_rates_2008}. We will use this result to establish the  the posterior convergence conditions for \gls{CGP}. The proof is deferred to Section \ref{sec:cgp_shift_ball_proof}.

Using Lemma \ref{thm:feasibility_ineq}, we are ready to establish the posterior convergence conditions for \gls{CGP}:

\begin{proof}
Show Theorem \ref{thm:cgp_conv} by showing conditions (\ref{eq:conv_cond_1})-(\ref{eq:conv_cond_3}) are satisfied for suitable choices of $\epsilon_n$ and $B_n$.
\begin{itemize}
\item Condition (\ref{eq:conv_cond_1}) is a direct consequence of Lemma \ref{thm:feasibility_ineq} and the assumption on concentration function. 

Recall that the definition for concentration function is $\psi_{g^*}(\epsilon) = \inf_{\hat{g} \in \Hsc_k, ||\hat{g} - g^*||_\infty \leq \epsilon} ||\hat{g}||^2_{\Hsc_k} - log \, P(||g||_\infty \leq \epsilon)$. Specifically, select $\hat{g}^*$ such that $||\hat{g}^* - g^*||<\epsilon_n$ and $||\hat{g}^*||^2_{\Hsc_k} =\inf_{\hat{g} \in \Hsc_k, ||\hat{g} - g^*||_\infty \leq \epsilon} ||\hat{g}||^2_{\Hsc_k}$, then $||g - g^*||_\infty \leq \epsilon_n + ||g - \hat{g}^*||_\infty$ and hence
\begin{align*}
P(||g - g^*||_\infty \leq 2\epsilon_n) 
&\geq P(||g - \hat{g}^*||_\infty \leq \epsilon_n)\\
&\geq exp(-\frac{1}{2}||\hat{g}^*||^2_{\Hsc_k}) * P(||g||_\infty \leq \epsilon_n)\\
&= exp(-\psi_{g^*}(\epsilon_n)) \\
& \geq exp(-n\epsilon_n^2)
\end{align*}
where the second inequality follows by Lemma \ref{thm:feasibility_ineq}, and the third inequality follows by assumption on concentration function $\psi(\epsilon_n) \leq n\epsilon_n^2$.

\item Conditions (\ref{eq:conv_cond_2}) and (\ref{eq:conv_cond_3}) can be shown by construction, i.e., by showing that there exists a set $B_n$ which satisfies  (\ref{eq:conv_cond_2}) and (\ref{eq:conv_cond_3}).\\

Denote $\Hbb_1$ a unit ball in the \gls{RKHS} $\Hsc_k$, denote $\Cbb_n$ the set of $\epsilon_n$-feasible functions such that $f(\Csc|g) > \alpha_n$, where $\alpha_n$ is a small positive constant, and denote $M_n$ a large postive constant. 

We construct $B_n$ as:
\begin{align*}
B_n &= M_n \Hbb_1 \cap \Cbb_n,
\end{align*}
that is, $B_n$ is a subset of the \gls{RKHS} that contains highly feasible functions. We show that $B_n$ satisfies (\ref{eq:conv_cond_2}) and (\ref{eq:conv_cond_3}) for suitable choices of $\alpha_n$ and $M_n$.

\begin{itemize}
\item First show $B_n$ satisfies (\ref{eq:conv_cond_2}). 

Notice that by De Morgan's law, $(M_n \Hbb_1 \cap \Cbb_n)^c = (M_n \Hbb_1)^c \cup \Cbb_n^c$, which implies that $P(g \not\in B_n) \leq P(g \not\in M_n \Hbb_1) + P(g \not\in \Cbb_n)$, therefore only need to show $P(g \not\in M_n \Hbb_1) + P(g \not\in \Cbb_n) \leq e^{-n\epsilon_n^2}$ for suitable choice of $M_n$ and $\alpha_n$.

First consider $P(g \not\in M_n \Hbb_1)$:
\begin{align*}
P(g \not\in M_n \Hbb_1) &= 
\int_{||g||_{\Hsc_k} > M_n} 
\frac{\phi_k(g)f(\Csc|g)}{E\big( f(\Csc|g) \big)}dg \\
& \leq
\int_{||g||_{\Hsc_k} > M_n} 
\frac{\phi_k(g)}{E\big( f(\Csc|g) \big)}dg \\
&= D * \int_{||g||_{\Hsc_k} > M_n} \phi_k(g) dg 
= D * P(||g||_{\Hsc_k} > M_n) \\
& \leq D * E\Big(exp(||g||^2_{\Hsc_k})\Big) * exp(-M_n^2) \\
& \leq D * G * exp(-M_n^2)
\end{align*}
where we have denoted $D, G$ two positive constants such that $D = \frac{1}{E\big( f(\Csc|g) \big)}$ and $G=E\Big(exp(||g||^2_{\Hsc_k})\Big)$. In the above equation, the first inequality follows since $f(\Csc|g) \leq 1$, and second inequality follows by first square and exponentiate both sides in $P(||g||_{\Hsc_k} > M_n)$, and then apply the Markov's inequality.

Now consider $P(g \not\in \Cbb_n)$:
\begin{align*}
P(g \not\in \Cbb_n) &= 
\int_{f(\Csc|g) < \alpha_n} 
\frac{\phi_k(g)f(\Csc|g)}{E\big( f(\Csc|g) \big)}dg \\
&\leq 
\frac{\alpha_n}{E\big( f(\Csc|g) \big)}
\int_{f(\Csc|g) < \alpha_n} \phi_k(g) dg\\
& =
\frac{\alpha_n}{E\big( f(\Csc|g) \big)} \\
& = D*\alpha_n
\end{align*}
where the first inequality follows by the definition of $\Cbb_n$, the second inequality follows since $\int_{f(\Csc|g) < \alpha_n} \phi_k(g) dg \leq 1$, i.e. integrating a Gaussian measure $\phi_k(g)$ over a subset of its full support yields a value less than 1.

Now for any $C>1$ such that $exp(-Cn\epsilon_n^2) \leq \frac{1}{2} * \frac{1}{D} * \frac{1}{G \vee 1} = \frac{1}{2} \frac{E\big( f(\Csc|g) \big)}{E\big(exp(||g||^2_{\Hsc_k})\big) \vee 1 }$, set
\begin{align}
M_n^2 = 2Cn \epsilon_n^2, \alpha_n = e^{-2Cn\epsilon_n^2},
\label{eq:conv_cond_m_a_val}
\end{align}
we then have 
\begin{align*}
P(g \not\in B_n) &\leq P(g \not\in M_n\Hbb_1) + P(g \not\in \Cbb_n) \\
&\leq D*G* exp(-M_n^2) + D * \alpha_n \\
&\leq \frac{1}{2} e^{-Cn\epsilon_n^2} + \frac{1}{2} e^{-Cn\epsilon_n^2}
= e^{-Cn\epsilon_n^2}.
\end{align*}
Therefore condition (\ref{eq:conv_cond_2}) is satisfied.

\item Now show $B_n$ satisfies the entropy number condition (\ref{eq:conv_cond_3}).

Recall that for the set $B_n$, its entropy number  $N(\epsilon_n, B_n, ||.||_\infty)$ is defined as the minimum number of balls of radius $\epsilon_n$ needed to cover $B_n$ in a metric space with norm $||.||_\infty$. Define $h_1, \dots, h_N$ elements of $M_n\Hbb_1 \cap \Cbb_n$ that are $2\epsilon_n$-separated with respect to the uniform norm $||.||_\infty$, then $||.||_\infty$-balls with radius $\epsilon_n$ and centers at $h_j$  are mutually disjoint. Therefore, denote $E_n = \{g|||g||_\infty < \epsilon_n\}$, for $g \in B_n$, we have:
\begin{align*}
1 
&\geq \sum_{j=1}^N P(g \in (h_j + E_n) \cap \Cbb_n) \\
&\geq \sum_{j=1}^N exp\Big(-\frac{1}{2}||h_j||^2_{\Hsc_k} \Big) 
P(g \in E_n)\\
&\geq N * exp\Big( -\frac{1}{2}M_n^2 \Big) 
P(g \in E_n)
\end{align*}
where the second inequality follows by Lemma \ref{thm:feasibility_ineq}, the third inequality follows since $h_j \in M_n \Hbb_1$ which implies $||h_j||^2_{\Hsc_k} \leq M_n$. It then follows that:
\begin{align*}
log N 
&\leq 
\frac{1}{2}M_n^2 - log \, P(g \in E_n) \\
& \leq C n \epsilon_n^2 - log \, P(g \in E_n)\\
& \leq C n \epsilon_n^2 + n \epsilon_n^2 \\
& \leq 2 Cn\epsilon_n^2,
\end{align*}
where the first inequality follows by the definition of $M_n$ in (\ref{eq:conv_cond_m_a_val}), second inequality follows by the assumption that the concentration function satisfy $\psi(\epsilon_n) \leq n\epsilon_n^2$, the last inequality follows since $C>1$.

Finally, let $h_1, \dots, h_N$ be maximal in the set $M_n \Hbb_1$ and recall that $h_j$'s are $2\epsilon_n$ separated, then by the definition of entropy number, we have:
\begin{align*}
log N(2\epsilon_n, B_n, ||.||_\infty) \leq log N 
&\leq 3 Cn\epsilon_n^2.
\end{align*}
Therefore condition (\ref{eq:conv_cond_3}) is satisfied.
\end{itemize}

\end{itemize}

\end{proof}

\subsection{Proof for Lemma \ref{thm:feasibility_ineq}}
\label{sec:cgp_shift_ball_proof}
\begin{proof}
Notice that the statement in the lemma is equivalent to:
\begin{align}
P(S_\epsilon + g^*) \geq exp\Big(-\frac{1}{2}||g^*||^2_{\Hsc_k} \Big) * P(S_\epsilon)
\label{eq:feas_state}
\end{align}
therefore only need to show above statement is true. 

\begin{itemize}
\item First show (\ref{eq:feas_state}) holds for $g^*$ that is $\epsilon$-feasible almost surely. 

Recall that almost surely $\epsilon$ feasibility is defined as $f(\Csc|g^* + s) \geq E(f(\Csc|s)) \; \forall s \in S_\epsilon$. Start from the left hand side of (\ref{eq:feas_state}), by the definition of zero-mean \gls{CGP}:
\begin{align*}
P(S_\epsilon + g^*) 
&= \int_{g' \in S_\epsilon + g^*} P(g') dg'
= \int_{g \in S_\epsilon} P(g + g^*) dg \\
&= \int_{g \in S_\epsilon} 
\frac{\phi_k(g + g^*)f(\Csc|g + g^*)}{E(f(\Csc|g))} dg\\
&\geq 
\frac{E(f(\Csc|s))}{E(f(\Csc|g))} 
*
\int_{g \in S_\epsilon} \phi_k(g + g^*) * dg \\
&=
\int_{z \in S_\epsilon} \frac{\phi_k(z)f(\Csc|z)}{\phi_k(S_\epsilon)E(f(\Csc|g))}
*
\int_{g \in S_\epsilon} \phi_k(g + g^*) * dg \\
&=
\frac{1}{\phi_k(S_\epsilon)} *
\int_{z \in S_\epsilon} \frac{\phi_k(z)f(\Csc|z)}{E(f(\Csc|g))}
*
\int_{g \in S_\epsilon} \phi_k(g + g^*) * dg \\
&=
\frac{P(S_\epsilon)}{\phi_k(S_\epsilon)} 
*
\int_{g \in S_\epsilon} \phi_k(g + g^*) * dg \\
&\geq 
P(S_\epsilon) * exp\Big(-\frac{1}{2}||g^*||^2_{\Hsc_k}\Big) 
\end{align*}
where the first equality follows from change of variables, the second equality follows from the definition of \gls{CGP}. The first inequality follows from the fact that $g^*$ is strongly $\epsilon$-feasible, the second inequality follows since $E(f(\Csc|g)) \leq 1$, the last inequality follows by the shift-ball inequality for Gaussian measures (Theorem 2 of \cite{kuelbs_gaussian_1994}). \\

\item Now show (\ref{eq:feas_state}) also holds for $g^*$ that is $\epsilon$-feasible in probability. 

Denote the event that a function $g$ violates the strong $\epsilon$ feasibility as $R_g = \Big\{g \Big| f(\Csc|g) \leq  \alpha_\epsilon \Big\}$, and the event that a noise $s$ makes $g^*$ violating the strong $\epsilon$ feasibility as $R_\epsilon = \Big\{\epsilon \Big| f(\Csc|g^* + \epsilon) \leq  \alpha_\epsilon \Big\}$. Also denote $g^* + s \sim truncGP(g^*, k, g^* + S_\epsilon)$ and $g^* + s' \sim CGP(g^*, k, \Csc)$ two perturbed functions by random noises distributed as truncated and constrained \gls{GP}s, respectively. Recall that weak $\epsilon$ feasibility upper bounds the probability of $g^* + s$ violating feasibility by requiring below to be true:
\begin{align*}
P(g^* + s \in R_g) 
\leq 
P(g^* + s' \in R_{g} | s' \in S_\epsilon)
\end{align*} 
First derive a useful fact using the definition of weak $\epsilon$ feasibility. Notice that:
\begin{align*}
P(g^* + s \in R_g) 
&= 
\int_{g^* + s \in R_g} \frac{\phi(g^* + s)}{\phi(g^* + S_\epsilon)} d(g^* + s)
\\
&= 
\int_{s \in R_\epsilon} \frac{\phi(g^* + s)}{\phi(g^* + S_\epsilon)} ds
\\
P(g^* + s' \in R_{g} | s' \in S_\epsilon)
&=
\frac{P(g^* + s' \in R_g, \quad g^* + s' \in g^* + S_\epsilon)}
{P(s' \in S_\epsilon)}
\\
&= 
\frac{1}{P(s' \in  S_\epsilon)} 
\int_{g^* + s' \in R_g \cap g^* + S_\epsilon} 
\frac{\phi(g^* + s')f(\Csc|g^* + s')}{E(f(\Csc|g))}d(g^* + s')
\\
&=
\frac{1}{P(s' \in  S_\epsilon)} 
\int_{s' \in R_\epsilon \cap S_\epsilon}
\frac{\phi(g^* + s')f(\Csc|g^* + s')}{E(f(\Csc|g))}d s'
\end{align*}
we can then derive below fact:
\begin{alignat}{3}
\int_{s' \in R_\epsilon \cap S_\epsilon}
\frac{\phi(g^* + s')f(\Csc|g^* + s')}{E(f(\Csc|g))}d s'
\geq 
P(s' \in  S_\epsilon) 
\int_{s \in R_\epsilon} \frac{\phi(g^* + s)}{\phi(g^* + S_\epsilon)} ds.
\label{eq:weak_feas_cond}
\end{alignat}

Using above fact, we are now ready to show (\ref{eq:feas_state}). Starting from the left-hand side:
\begin{align}
P(S_\epsilon + g^*) 
&= \int_{g' \in S_\epsilon + g^*} P(g') dg'
= \int_{g \in S_\epsilon} P(g + g^*) dg 
\nonumber\\
&= 
\int_{g \in S_\epsilon \cap R_\epsilon^c} P(g + g^*) dg + 
\int_{g \in S_\epsilon \cap R_\epsilon} P(g + g^*) dg
\nonumber\\
&= 
\int_{g \in S_\epsilon \cap R_\epsilon} 
\frac{\phi_k(g + g^*)f(\Csc|g + g^*)}{E(f(\Csc|g))} dg +
\int_{g \in S_\epsilon \cap R_\epsilon^c} 
\frac{\phi_k(g + g^*)f(\Csc|g + g^*)}{E(f(\Csc|g))} dg 
\nonumber\\
& \geq 
\int_{g \in S_\epsilon \cap R_\epsilon} 
\frac{\phi_k(g + g^*)f(\Csc|g + g^*)}{E(f(\Csc|g))} dg +
\frac{P(S_\epsilon)}{\phi_k(S_\epsilon)}
\int_{g \in S_\epsilon \cap R_\epsilon^c} \phi_k(g + g^*) dg
\label{eq:weak_feas_midpoint}
\end{align}
where the last inequality follows by noticing $S_\epsilon \cap R_\epsilon^c$ is the region where the strong $\epsilon$-feasibility $f(\Csc|g + g^*) \geq E(f(\Csc|s))$ is satisfied, therefore the inequality in the second integral follows by the definition of strong $\epsilon$-feasibility (see the proof for strong $\epsilon$-feasibility for detail). We now handle the first integral by applying the fact (\ref{eq:weak_feas_cond}):
\begin{align*}
(\ref{eq:weak_feas_midpoint}) 
& =
\int_{g \in S_\epsilon \cap R_\epsilon} 
\frac{\phi_k(g + g^*)f(\Csc|g + g^*)}{E(f(\Csc|g))} dg +
\frac{P(S_\epsilon)}{\phi_k(S_\epsilon)}
\int_{g \in S_\epsilon \cap R_\epsilon^c} \phi_k(g + g^*) dg
\\
& \geq 
\frac{P(S_\epsilon)}{\phi_k(S_\epsilon + g^*)}
\int_{g \in S_\epsilon \cap R_\epsilon} \phi_k(g + g^*) dg + 
\frac{P(S_\epsilon)}{\phi_k(S_\epsilon)}
\int_{g \in S_\epsilon \cap R_\epsilon^c} \phi_k(g + g^*) dg 
\\
& \geq 
\frac{P(S_\epsilon)}{\phi_k(S_\epsilon)}
\int_{g \in S_\epsilon \cap R_\epsilon} \phi_k(g + g^*) dg + 
\frac{P(S_\epsilon)}{\phi_k(S_\epsilon)}
\int_{g \in S_\epsilon \cap R_\epsilon^c} \phi_k(g + g^*) dg 
\\
& = 
\frac{P(S_\epsilon)}{\phi_k(S_\epsilon)} 
\int_{g \in S_\epsilon} \phi_k(g + g^*) dg \\
& \geq 
P(S_\epsilon) * exp\Big(-\frac{1}{2}||g^*||^2_{\Hsc_k}\Big),
\end{align*}
where the first inequality follows by the fact (\ref{eq:weak_feas_cond}), the second inequality follows by the Anderson's theorem \citep{anderson_integral_1955, gardner_brunn-minkowski_2002}, and the last inequality follows by the shift-ball inequality for Gaussian measures (Theorem 2 of \cite{kuelbs_gaussian_1994}).

\item Finally show (\ref{eq:feas_state}) holds for $g^*$ that is $\epsilon$-feasible in expectation. 

Recall that the in-expectation $\epsilon$ feasibility is defined as $E(f(\Csc|g^* + s)) \geq E(f(\Csc|s'))$ for $g^* + s \sim truncGP(g^*, k, g^*+S_\epsilon)$ and $s' \sim truncGP(0, k, S_\epsilon)$. Also notice that:
\begin{align*}
E(f(\Csc|g^* + s)) &= 
\int_{s \in S_\epsilon} f(\Csc|g^* + s) \frac{\phi(g^* + s)}{\phi(g^* + S_\epsilon)} ds\\
E(f(\Csc|s)) &=
\int_{s \in S_\epsilon} f(\Csc|s) \frac{\phi(s)}{\phi(S_\epsilon)} ds,
\end{align*}
therefore the in-expectation $\epsilon$ feasibility implies that
\begin{align}
\int_{s \in S_\epsilon} f(\Csc|g^* + s) \phi(g^* + s) ds
&  \geq 
\frac{\phi(g^* + S_\epsilon)}{\phi(S_\epsilon)} * \int_{s \in S_\epsilon} f(\Csc|s) \phi(s) ds
\label{eq:in_expect_feas_fact}
\end{align}
We are now ready to show (\ref{eq:feas_state}), start from the left hand side:
\begin{align*}
P(S_\epsilon + g^*) 
&= \int_{g' \in S_\epsilon + g^*} P(g') dg'
= \int_{g \in S_\epsilon} P(g + g^*) dg \\
&= \int_{g \in S_\epsilon} 
\frac{\phi_k(g + g^*)f(\Csc|g + g^*)}{E(f(\Csc|g))} dg\\
&\geq 
\frac{\phi(g^* + S_\epsilon)}{\phi(S_\epsilon)} 
* 
\int_{s \in S_\epsilon} \frac{f(\Csc|s) \phi(s)}{E(f(\Csc|g))} ds\\
&=
\frac{\phi(g^* + S_\epsilon)}{\phi(S_\epsilon)} * P(S_\epsilon) \\
&\geq 
exp\Big(-\frac{1}{2}||g^*||^2_{\Hsc_k}\Big) * P(S_\epsilon)
\end{align*}
where the first inequality follows from the fact (\ref{eq:in_expect_feas_fact}) implied by the in-expectation $\epsilon$-feasiblity, the second inequality follows by the shift-ball inequality for Gaussian measures (Theorem 2 of \cite{kuelbs_gaussian_1994}).

\end{itemize}

\end{proof}

\end{document}